\documentclass[12pt]{amsart}
\usepackage{graphicx,color}
\usepackage{colordvi}
\usepackage{amssymb, mathrsfs, amsfonts, amsmath}
\usepackage{latexsym}

\newtheorem{theorem}{Theorem}[section]
\newtheorem{lemma}[theorem]{Lemma}

\newtheorem{proposition}[theorem]{Proposition}
\newtheorem{definition}[theorem]{Definition}

\newtheorem{remark}[theorem]{Remark}

\newcommand{\C}{\mathbb C}
\newcommand{\R}{\mathbb R}
\newcommand{\Q}{\mathbb Q}

\begin{document}

\title[Lipschitz Geometry of Surfaces
 ]{Lipschitz Geometry of Real Semialgebraic Surfaces}

\author[]{Lev Birbrair}
\address{Departamento de Matem\'atica, Universidade Federal do Cear\'a
(UFC), Campus do Pici, Bloco 914, Cep. 60455-760. Fortaleza-Ce,
Brasil and \newline Institute of Mathematics, Jagiellonian University, Profesora Stanisława Łojasiewicza 6, 30-348 Kraków, Poland  }
 \email{lev.birbrair@gmail.com}
\author[]{Andrei Gabrielov}
\address{Department of Mathematics, Purdue University,
	West Lafayette, IN 47907, USA}\email{gabrielov@purdue.edu}


\begin{abstract}
We present here basic results in Lipschitz Geometry of semialgebraic surface germs. Although bi-Lipschitz classification problem of surface germs with respect to the inner metric was solved long ago, classification with respect to the outer metric remains an open problem. We review recent results related to the outer and ambient bi-Lipschitz classification of surface germs. In particular, we explain why the outer Lipschitz classification is much harder than the inner classification, and why the ambient Lipschitz Geometry of surface germs is very different from their outer Lipschitz Geometry. In particular, we show that the ambient Lipschitz Geometry of surface germs includes all of the Knot Theory.
\end{abstract}

\maketitle

\section{Introduction}
Lipschitz classification of semialgebraic surfaces has become in recent years one of the central questions of the
Metric Geometry of Singularities.
It was stimulated by the finiteness theorems of Mostowski, Parusinski and Valette.
(see \cite{Mostowski, Parusinski, valette2005Lip}).
They proved that there are finitely many Lipschitz equivalence classes in any semialgebraic family of semialgebraic sets. Lipschitz classification is intermediate between Smooth (too fine) and Topological (too coarse) classifications.
For example, smooth classification of most singularities is not finite.
It may be even infinite dimensional for non-isolated singularities.

Here we review recent developments in Lipschitz Geometry of semialgebraic surfaces (two-dimensional real semialgebraic sets).
Since we are mainly interested in singularities of semialgebraic surfaces, our main object is a
semialgebraic surface germ $(X,0)$ at the origin of $\R^n$. Note that most results presented in this paper
remain true for subanalytic sets, and for the sets definable in a polynomially bounded o-minimal structure.

A connected semialgebraic set $X\subset\R^n$ inherits from $\R^n$ two metrics:
the {\bf outer metric} $dist(x,y) = |y-x|$ and the {\bf inner metric}
$idist(x,y) =$ length of the shortest path in $X$ connecting $x$ and $y$.
Note that $dist(x,y)\le idist(x,y)$.
A semialgebraic set is called Lipschitz Normally Embedded if these two metrics are equivalent (see Definition \ref{LNE}).

For the surface germs, there are three natural equivalence relations:

{\bf 1)} Inner Lipschitz equivalence: $(X,0)\sim_i (Y,0)$ if there is a homeomorphism\newline
$h:(X,0)\to(Y,0)$ bi-Lipschitz with respect to the inner metric.

{\bf 2)} Outer Lipschitz equivalence: $(X,0)\sim_o (Y,0)$ if there is a homeomorphism\newline
$h:(X,0)\to(Y,0)$ bi-Lipschitz with respect to the outer metric.

{\bf 3)} Ambient Lipschitz equivalence:
$(X,0)\sim_a (Y,0)$ if there is an orientation preserving bi-Lipschitz homeomorphism $H:(\R^n,0)\to(\R^n,0)$ such that $H(X)=Y$.

Inner Lipschitz Geometry of surface germs is relatively simple.
The building block of the inner Lipschitz classification of surface germs is a $\beta$-{\bf H\"older triangle} (see Definition \ref{holder}).
A surface germ $(X,0)$ with an isolated singularity is inner Lipschitz equivalent to a $\beta$-{\bf horn} (see Definition \ref{horn}).
If the singularity is not isolated, classification is made by the theory of H\"older Complexes (see \cite{Birbrair}).
A {\bf H\"older Complex} is a triangulation (decomposition into H\"older triangles) of a surface germ.
Canonical H\"older Complex (see Definition \ref{canonicalHC}) is a complete invariant of the inner Lipschitz equivalence of surface germs.

For example, the germs of all irreducible complex curves are inner Lipschitz equivalent to $(\C,0)$, while the outer Lipschitz classification of the germs of complex plane curves is described by their sets of essential Puiseux pairs (see \cite{Pham-Teissier}, \cite{fernandes}). Even for the union of two normally embedded H\"older triangles, the outer Lipschitz Geometry is not simple (see \cite{BG}).

A special case of a surface germ is the union of a H\"older triangle $T$ and a graph of a Lipschitz semialgebraic function $f$ defined on $T$. The outer Lipschitz equivalence of two such surface germs is equivalent to the Lipschitz contact equivalence of the two functions. This relates outer Lipschitz geometry of surface germs with the Lipschitz geometry of functions. In \cite{birbrair2014lipschitz} a complete invariant of the contact equivalence class of a Lipschitz function $f$ defined on a H\"older triangle $T$, called a ``pizza,'' is defined. Informally, a pizza is a decomposition of $T$ into ``slices,'' H\"older sub-triangles  $\{T_i\}$ of $T$, such that the order of $f$ on each arc $\gamma\subset T_i$ depends linearly on the order of contact of $\gamma$ with a boundary arc of $T_i$.

For the general surface germs, Lipschitz classification with respect to the outer metric is still an open problem. The set of semialgebraic arcs in $(X,0)$ parameterized by the distance to $0$ is called the Valette link $V(X)$ of the germ $(X,0)$ (see Definition \ref{arc}). The order of contact of the arcs (see Definition \ref{tord}) defines a non-archimedean metric on $V(X)$.

The study of Lipschitz Normally Embedded, or simply Normally Embedded, sets was initiated by Kurdyka and Orro \cite{kurdyka-orro}.
Kurdyka proved that any semialgebraic set admits a finite partition into normally embedded subsets.
Using this partition, Kurdyka and Orro proved that any semialgebraic sets admits a semialgebraic ``pancake metric'' equivalent to the inner metric. Normal Embedding theorem of Birbrair and Mostowski states that, for any semialgebraic set $X$, there is a semialgebraic and bi-Lipschitz with respect to the inner metric embedding $\Psi: X \to \R^m$, where $m\ge 2\dim(X)+1$ (see \cite{extention}).
Lipschitz Normal Embedding of complex analytic sets is addressed in the paper by Anne Pichon in the present volume.

A pancake decomposition is called minimal if it is not a refinement of another pancake decomposition.
A natural question related to Lipschitz Normal Embedding of surface germs is uniqueness of a minimal pancake decomposition.
The answer is negative even for a H\"older triangle.
Gabrielov and Sousa in \cite{GS} gave examples of H\"older triangles having several combinatorially non-equivalent minimal pancake decompositions.

Relations between ambient and outer equivalence of surface germs were studied in \cite{met-knots1,met-knots2,met-knots3}. In the paper  \cite{met-knots1} the authors presented several outer Lipschitz and ambient topologically equivalent families of surface germs $(X_i,0)$, which were pairwise ambient Lipschitz non-equivalent. In \cite{met-knots2,met-knots3},  several ``Universality Theorems'' were formulated.
Informally, these theorems state that, even when the link of a surface germ is topologically a trivial knot, ambient Lipschitz classification of such surface germs ``contains all of the Knot Theory.''

\section{Inner Lipschitz Equivalence}\label{sect:inner}

\begin{figure}
\centering
\includegraphics[width=6in]{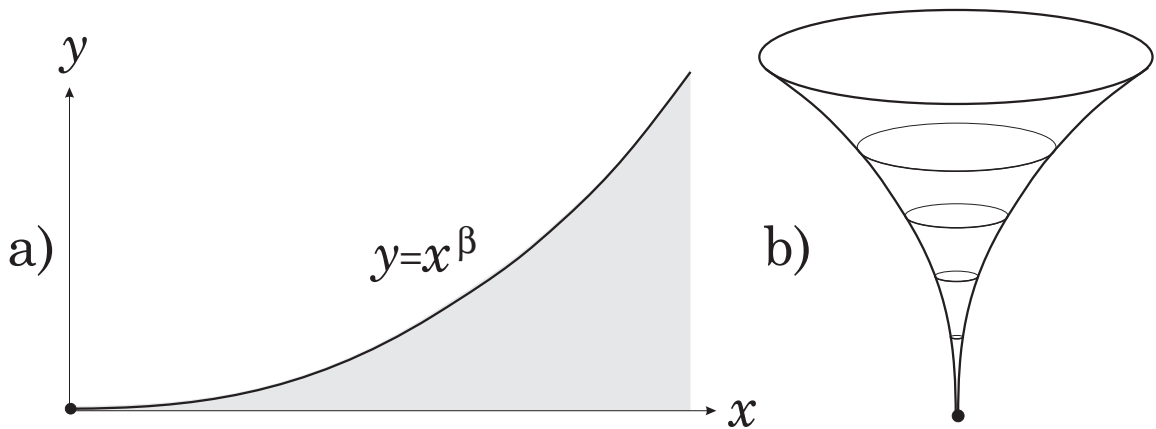}
\caption{A $\beta$-H\"older triangle and a $\beta$-horn.}\label{holder-horn}
\end{figure}

\begin{definition}\label{holder}\normalfont
For $1\le\beta \in \Q$, the standard $\beta$-H\"older triangle $T_\beta$ is the germ at the origin of $\R^2$ of the surface
 $\{x\ge 0, \; 0\le y \le x^\beta\}$ (see Figure \ref{holder-horn}a). A $\beta$-\emph{H\"older triangle} is a surface germ inner Lipschitz equivalent to $T_\beta$.
\end{definition}

\begin{definition}\label{horn}\normalfont
For $1\le\beta \in \Q$, the standard $\beta$-horn $C_\beta$ is the germ at the origin of $\R^3$ of the surface $\{z\ge 0, \; x^2+y^2 = z^{2\beta}\}$ (see Figure \ref{holder-horn}b).
 A $\beta$-\emph{horn} is a surface germ inner Lipschitz equivalent to $C_\beta$.
\end{definition}

 \begin{theorem}\label{thm:Horn-theorem}
	Given the germ $(X,0)$ of a semialgebraic surface with isolated singularity and connected link, there is a unique rational
number $\beta\geq 1$ such that $(X,0)$ is inner Lipschitz equivalent to the standard $\beta$-horn $C_\beta$.
\end{theorem}

Birbrair's theory of H\"older Complexes (see \cite{Birbrair}) is a generalization of Theorem
\ref{thm:Horn-theorem} for the surface germs with non-isolated singularities.

\begin{definition}\label{Hol-com} \normalfont A \emph{Formal H\"older Complex} is a pair $(G,\beta)$, where $G$ is a graph
 and $\beta: E_G\to\Q_{\ge 1}$ is a function, where $E_G$ the set of edges of $G$.
\end{definition}

\begin{definition}\label{geometric-HC}\normalfont A \emph{Geometric H\"older Complex} corresponding to a Formal
H\"older Complex $(G,\beta)$ is a surface germ $(X,0)$ such that\newline
{\bf 1.} For small $\varepsilon>0$, the intersection of $X$ with the $\varepsilon$-ball $B_\varepsilon$ is homeomorphic to the cone over $G$, and the intersection of $X$ with the $\varepsilon$-sphere $S_\varepsilon$ is homeomorphic to $G$.\newline
{\bf 2.} For any edge $g\in E_G$, the subgerm of $(X,0)$ corresponding to $g$ is a $\beta(g)$-H\"older triangle.
\end{definition}

\begin{figure}
\centering
\includegraphics[width=6in]{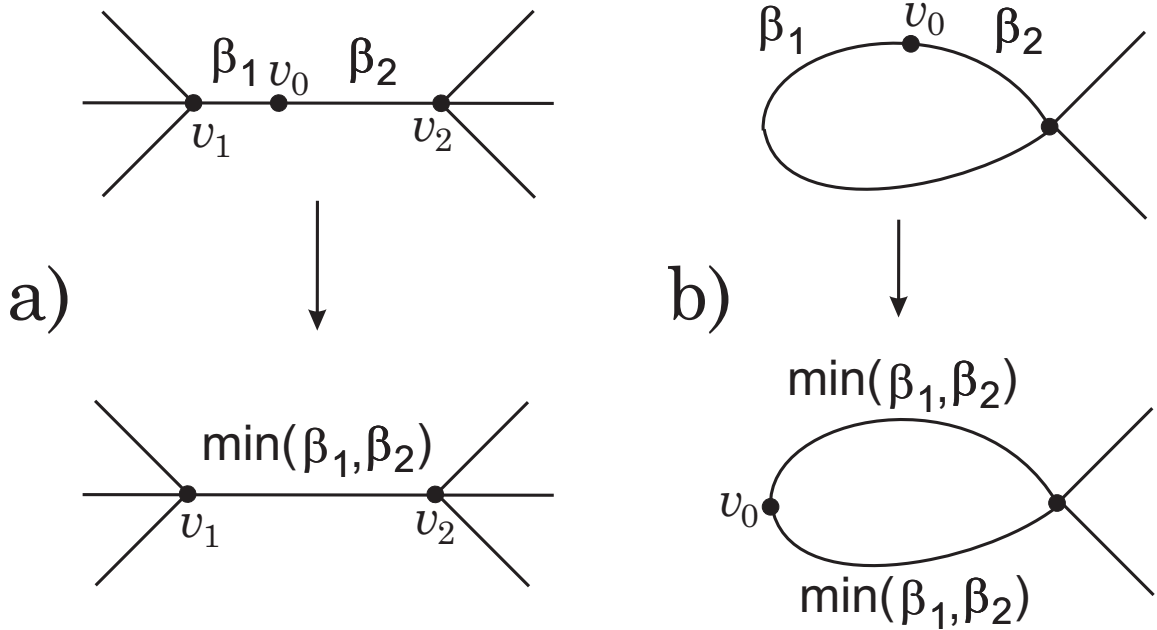}
\caption{Simplification of H\"older complexes}\label{simplification-yes}
\end{figure}

\begin{theorem}\label{surf-Hold}
For any surface germ $(X,0)\subset \R^n$, there exists a Formal H\"older Complex $(G,\beta)$ such that
$(X,0)$ is a Geometric H\"older Complex corresponding to $(G,\beta)$.
\end{theorem}

\begin{remark}\label{needs-simplification} \normalfont For a given surface germ $(X,0)$, the Formal H\"older complex $(G,\beta)$ in Theorem \ref{surf-Hold} is not unique. The simplification procedure described below reduces it to the unique Canonical H\"older Complex corresponding to $(X,0)$.
\end{remark}

\begin{definition}\label{simplification} \normalfont We say that a vertex $v_0$ of the graph $G$ is {\it non-critical} if
it is adjacent to exactly two edges $g_1$ and $g_2$ of $G$, and these edges connect $v_0$ with two different vertices of $G$.
A vertex $v_0$ of $G$ is called a {\it loop vertex} if it is adjacent to exactly two different edges $g_1$ and $g_2$ of $G$,
and these edges connect $v_0$ with the same vertex $v_1$ of $G$.
The other vertices of $G$ (neither non-critical nor loop vertices) are called {\it critical}.
\end{definition}

\begin{definition}\label{canonicalHC} \normalfont An Abstract H\"older Complex $(G,\beta)$ is called \emph{Canonical}, if\newline
{\bf 1.} All  vertices of $G$ are either critical or loop vertices;\newline
{\bf 2.} For any loop vertex $v$ of $G$ adjacent to the edges $g_1$ and $g_2$, one has $\beta(g_1)=\beta(g_2)$.
\end{definition}

Now we define a {\bf simplification procedure}, reducing an Abstract H\"older Complex $(G,\beta)$ to a Canonical one.

We start with eliminating non-critical vertices.
Let $v_0$ be a non-critical vertex of $G$, connected with the vertices $v_1$ and $v_2$ of $G$ by the edges $g_1$ and $g_2$.
Then we remove the vertex $v_0$ from $V(G)$, and replace the edges $g_1$ and $g_2$ of $G$ with the single edge $g_0$ connecting  $v_1$ with $v_2$.
Let $G'$ be the graph obtained from $G$ by this operation.
We define an abstract H\"older complex $(G',\beta')$, setting $\beta'(g_0)=\min\{\beta(g_1),\beta(g_2)\}$ and
$\beta'(g)=\beta(g)$ on all edges $g$ of $G'$ other than $g_0$, which are also the edges of $G$. (see Figure \ref{simplification-yes}a)

We repeat this operation until there are no non-critical vertices. After that, we take care of the loop vertices of $G$.

Let $(G,\beta)$ be an Abstract H\"older Complex without non-critical vertices.
If a loop vertex $v_0$ of $G$ is connected by the edges $g_1$ and $g_2$ with the same vertex $v_1$,
such that $\beta(g_1)\ne\beta(g_2)$, we define an Abstract H\"older Complex $(G,\beta')$,
replacing $\beta_1=\beta(g_1)$ and $\beta_2=\beta(g_2)$ with $\beta'(g_1)=\beta'(g_2)=\min(\beta_1,\beta_2)$
(see Figure \ref{simplification-yes}b).
We repeat this operation for all loop vertices of $G$.

The main results of the paper \cite{Birbrair} are the following:

\begin{theorem}\label{inner-classification} \emph({\bf Inner Lipschitz Classification Theorem.}\emph)
The surface-germs $(X,0)$ and $(X',0)$ are Lipschitz equivalent with respect to the inner metric if, and only if, the corresponding
Canonical H\"older Complexes are combinatorially equivalent.
\end{theorem}

\begin{theorem}\label{REALIZATION-THEOREM} \emph({\bf Realization Theorem.}\emph)
Let $(G,\beta)$ be an Abstract H\"older Complex. Then there
exists a surface germ $(X,0)$ which
 is a Geometric H\"older Complex corresponding to $(G,\beta)$.
\end{theorem}

 \begin{remark}\label{complexcurveinner}\normalfont
The theory of H\"older Complexes implies that a germ $(X,0)$ of an irreducible complex curve, considered as a real surface germ,
is inner Lipschitz equivalent to the germ $(\C,0)$.
Otherwise $(X,0)$ is inner Lipschitz equivalent to the union of finitely many germs $(\C,0)$ pinched at the origin, corresponding to irreducible components of $(X,0)$.
\end{remark}

\section{Normal Embedding Theorem, Lipschitz Normally Embedded Sets}

\begin{definition}\label{LNE}\normalfont A semialgebraic set $X$ is called \emph{Lipschitz Normally Embedded} (LNE) if the inner and outer metrics on $X$ are equivalent: $dist(x,y)\le idist(x,y)\le C\,dist(x,y)$ for some constant $C>0$ and all $x,y\in X$.
\end{definition}

For example, the germ of an algebraic curve $\{x^3=y^2\}$ is not LNE, while the standard $\beta$-horn $C_\beta$ is LNE.
 A germ of an irreducible complex curve is LNE if, and only if, it is smooth.

There are many examples of not normally embedded surface germs. On the other hand, we have the following result:

\begin{theorem}\label{LNE-theorem} \emph({See \cite{birbrair2000normal}.}\emph) Let $X\subset\R^m$ be a connected
  semialgebraic set. Then there exist a normally embedded
 semialgebraic set $\tilde{X}\subset\R^q$ and an
inner  bi-Lipschitz   homeomorphism $p:  X\to \tilde{X}$. This map is called a normal embedding of $X$.
\end{theorem}

\begin{definition}\label{LNE2}\normalfont A subset $\tilde{X}\subset\R^m$ is called Lipschitz Normally Embedded if there exist a bi-Lipschitz
homeomorphism $\Psi: \tilde{X}_{inner}
\rightarrow \tilde{X}_{outer}$.
\end{definition}

Here $\tilde{X}_{inner}$ means the space $\tilde{X}$ equipped with the inner metric, and $\tilde{X}_{outer}$ means $\tilde{X}$ equipped
with the outer metric. The difference with Definition \ref{LNE}
is that in Definition \ref{LNE2} we do not a priori suppose that $\Psi$ is the identity map.

\begin{proposition}\label{equivalence}
The two definitions of Lipschitz Normally Embedding are equivalent.
\end{proposition}

 Pancake Decomposition of Kurdyka implies that there exists a decomposition of any semialgebraic set
$(X,0)$ into LNE semialgebraic subsets.

\begin{theorem}\label{pancake-decomposition}\emph({See \cite{kurdyka-orro}.}\emph)
There is a decomposition of any semialgebraic set $X$ into subsets $X_i$ such
that :

1. $X_i$ are semialgebraic LNE sets.

2. $\dim(X_i \cap X_j)< \min(\dim X_i,\dim X_j)$ for $i\ne j$.
\end{theorem}

\begin{remark}\label{kurdyka&orro}\normalfont
Using pancake decomposition, Kurdyka and Orro defined the so-called
\emph{pancake metric} (see \cite{kurdyka-orro}, \cite{birbrair2000normal}).
It is a semialgebraic metric equivalent to the inner metric.
\end{remark}

\begin{definition}\label{arc}\normalfont (See \cite{Valette-link}.)
An \emph{arc} in $\R^n$ is (a germ at the origin of) a mapping $\gamma:[0,\epsilon)\rightarrow \R^n$ such that $\gamma(0)=0$.
Unless otherwise specified, arcs are parameterized by the distance to the origin, i.e., $|\gamma(t)|=t$.
We usually identify an arc $\gamma$ with its image in $\R^n$.
The \emph{Valette link} of a surface germ $(X,0)$ is the set $V(X)$ of all arcs $\gamma\subset X$.
\end{definition}

\begin{theorem}\label{valette-link-th} \emph{(See \cite{Valette-link}.)}
Let $(X,0)$ and $(Y,0)$ be germs of semialgebraic sets in $\R^n$.
If these germs are semialgebraically (inner, outer or ambient) Lipschitz equivalent,
then there exists a bi-Lipschitz map $h:X\to Y$ (or $h:\R^n\to\R^n$ such that $h(X)=Y$ in the case of ambient equivalence) such that $h(X\cap S_\varepsilon)=Y\cap S_\varepsilon$ for small $\varepsilon>0$.
\end{theorem}

 \begin{definition}\label{ordonarc}
\emph{Let $f\not\equiv 0$ be (a germ at the origin of) a Lipschitz function defined on an arc $\gamma$.
The \emph{order} of $f$ on $\gamma$ is $q=ord_\gamma f\in\Q$ such that $f(\gamma(t))=c t^q+o(t^q)$ as $t\to 0$, where $c\ne 0$.
If $f\equiv 0$ on $\gamma$, we set $ord_\gamma f=\infty$.}
\end{definition}

\begin{definition}\label{tord}
\emph{ The \emph{tangency order} of arcs $\gamma$ and $\gamma'$ is  $tord(\gamma,\gamma')=ord_{\gamma}
|\gamma(t)-\gamma'(t)|$.
The tangency order of an arc $\gamma$ and a set of arcs $Z\subset V(X)$  is
$tord(\gamma,Z)=\sup_{\lambda\in Z} tord(\gamma,\lambda)$.
The tangency order of two subsets $Z$ and $Z'$ of $V(X)$ is $tord(Z,Z')=\sup_{\gamma\in Z} tord(\gamma,Z')$.
Similarly, $itord(\gamma,\gamma')$, $itord(\gamma,Z)$ and $itord(Z,Z')$ are the tangency orders with respect to the inner metric.
 }
\end{definition}

\begin{remark}\label{semialgebraic-curve}\normalfont(See \cite{birbrair-fernandes-curves}.)
If $(X,0)$ is a germ of a semialgebraic curve, i.e., $X=\cup \gamma_i$ is a finite family of semialgebraic arcs, then the outer Lipschitz Geometry of $(X,0)$ is totaly determined by the tangency orders $\{tord(\gamma_i,\gamma_j)\}$.
\end{remark}

 \begin{proposition}\label{arc-criterion}
A surface germ $(X,0)$ is LNE if, and only if, for any two arcs $\gamma_1,\gamma_2$ in $X$ one has
$tord(\gamma_1,\gamma_2)=itord(\gamma_1,\gamma_2)$.
\end{proposition}

\begin{proposition}\label{generic-projection}
Let $(X,0) \subset (\R^n,0)$ be a $\beta$-horn. The Grassmannian $G(n,2)$ can be considered as the space of all orthogonal projections $\rho:\R^n\to\R^2$. Then there exist an open semialgebraic subset $\widetilde{G} \subset G(n,2)$ such that for all $\rho \in \widetilde{G}$ one has $\beta=\min_{\{\gamma_1,\gamma_2\}\subset V(X)} tord(\rho(\gamma_1),\rho(\gamma_2))$.
\end{proposition}

The following proposition was proved first by Alexandre Fernandes \cite{fernandes}. A special case of this is the Arc Criterion of Normal Embedding \cite{com-Rodrigo}.

\begin{proposition}\label{important-proposition}
 Let $(X,0)$ and $(Y,0)$ be surface germs. A semialgebraic homeomorphism $\Phi:(X,0)\to(Y,0)$ preserving
 the distance to the origin is outer bi-Lipschitz if, and only if, for any two arcs $\gamma_1,\gamma_2 \in V(X)$ one has
 \begin{equation}\label{twoarcs}
 tord(\gamma_1,\gamma_2)=tord(\Phi(\gamma_1),\Phi(\gamma_2)).
 \end{equation}
\end{proposition}

A special case of Pancake Decomposition for surface germs can be stated as follows:

\begin{theorem}\label{pancake-decomposition-surface}
Let $(X,0)$ be a surface germ. Then there exists a decomposition of $(X,0)$ into the germs $(X_i,0)$ such
that

1. Each $(X_i,0)$ is a LNE $\beta_i$-H\"older triangle.

2. For $i\ne j$, the intersection $(X_i,0) \cap (X_j,0)$ is either the origin or a common boundary arc of $(X_i,0)$ and $(X_j,0)$.
\end{theorem}

\begin{definition}\label{minimal-pancake}\normalfont
A pancake decomposition of a surface germ is \emph{minimal} if the union of any two adjacent H\"older triangles $X_i$ and $X_j$ is not normally embedded.
Two pancake decompositions are  \emph{combinatorially equivalent} if they are combinatorially equivalent as H\"older Complexes.
\end{definition}

The answer to a natural question ``Are any two minimal pancake decompositions of the same surface germ
combinatorially equivalent?'' is negative (see Section \ref{section:snakes}).

\section{Pizza Decomposition of the Germ of a Semialgebraic Function}\label{section:pizza}

This section is related to the outer Lipschitz Geometry of a special kind of a surface germ: The union of a LNE H\"older triangle and the graph
of a semialgebraic Lipschitz function defined on it.

\begin{definition}\label{def:Qf}
	\emph{For a semialgebraic Lipschitz function $f$ defined on a $\beta$-H\"older triangle $T$, let
	\begin{equation}\label{eqn:Qf}
	Q_{f}(T)=\bigcup_{\gamma\in V(T)}ord_{\gamma} f.
	\end{equation}
	It was shown in \cite{birbrair2014lipschitz} that $Q_f(T)$ is either a point or a closed interval in $\Q\cup\{\infty\}$.}
\end{definition}

\begin{definition}\label{def: Elementary Holder triangle}
	\emph{A H\"older triangle $T$ is \emph{elementary} with respect to a Lipschitz function $f$ if, for any $q\in Q_f(T)$ and any two arcs
$\gamma$ and $\gamma'$ in $T$ such that $ord_{\gamma}f=ord_{\gamma'}f=q$, the order of $f$ is $q$ on any arc in the H\"older triangle
$T(\gamma,\gamma')\subseteq T$ bounded by the arcs $\gamma$ and $\gamma'$.}
\end{definition}

\begin{definition}\label{def:width function}\normalfont Let $T$ be a H\"older triangle and $f$ a Lipschitz function defined on $T$.
For each arc $\gamma \subset T$, the \emph{width} $\mu_T(\gamma,f)$ of the arc $\gamma$ with respect to $f$ is the infimum of exponents of H\"older triangles $T'\subset T$ containing $\gamma$ such that $Q_f(T')$ is a point.
For $q\in Q_f(T)$ let $\mu_{T,f}(q)$ be the set of exponents $\mu_T(\gamma,f)$, where $\gamma$ is any arc in $T$ such that $ord_{\gamma}f=q$.
It was shown in \cite{birbrair2014lipschitz} that, for each $q\in Q_f(T)$, the set $\mu_{T,f}(q)$ is finite.
This defines a multivalued \emph{width function} $\mu_{T,f}: Q_f(T)\to \Q\cup \{\infty\}$.
If $T$ is elementary with respect to $f$, then the function $\mu_{T,f}$ is single valued.
When $f$ is fixed, we write $\mu_T(\gamma)$ and $\mu_T$ instead of $\mu_T(\gamma,f)$ and $\mu_{T,f}$.
\end{definition}

\begin{definition}\label{def: pizza slice}\normalfont
Let $T$ be a H\"older triangle and $f$ a semialgebraic Lipschitz function defined on $T$.
We say that $T$ is a \emph{pizza slice} associated with
$f$ if it is elementary with respect to $f$ and, unless $Q_f(T)$ is a point, $\mu_{T,f}(q)=aq+b$ is an affine function on $Q_f(T)$.
If $T$ is a pizza slice such that $Q_f(T)$ is not a point, then the \emph{supporting arc}
$\tilde\gamma$ of $T$ with respect to $f$ is the boundary arc of $T$ such that $\mu_T(\tilde\gamma,f)=\max_{q\in Q_f(T)}\mu_{T,f}(q)$.
In that case, $\mu_T(\gamma,f)=tord(\gamma,\tilde\gamma)$ for any arc $\gamma\subset T$ such that $tord(\gamma,\tilde\gamma)\le \mu_T(\tilde\gamma,f)$.
\end{definition}

\begin{definition}\label{pizza-decomp}\emph{(See \cite[Definition 2.13]{birbrair2014lipschitz}.)}
\emph{Let $f$ be a non-negative semialgebraic Lipschitz function defined on a $\beta$-H\"older triangle $T=T(\gamma_1,\gamma_2)$
oriented from $\gamma_1$ to $\gamma_2$.
A \emph{pizza} on $T$ associated with $f$ is a decomposition $\{T_\ell\}_{\ell=1}^p$ of $T$ into $\beta_\ell$-H\"older triangles $T_\ell=T(\lambda_{\ell-1},\lambda_\ell)$ ordered according to the orientation of $T$, such that
$\lambda_0=\gamma_1$ and $\lambda_p=\gamma_2$ are the boundary arcs of $T$,
$T_\ell\cap T_{\ell+1}=\lambda_\ell$ for $0<\ell< p$, and
each triangle $T_\ell$ is a pizza slice associated with $f$.}

\emph{A pizza $\{T_\ell\}$ on $T$ is \emph{minimal} if $T_{\ell-1}\cup T_\ell$ is not a pizza slice for any $\ell>1$.}
\end{definition}

\begin{definition}\label{abstract_pizza} \emph{(See \cite[Definition 2.12]{birbrair2014lipschitz}.)} \emph{An \emph{abstract pizza} is a finite ordered sequence $\{q_\ell\}_{\ell=0}^p$, where $q_\ell\in\Q_{\ge 1}\cup\{\infty\}$, and a finite collection $\{\beta_\ell,Q_\ell,\mu_\ell\}_{\ell=1}^p$, where $\beta_\ell\in\Q_{\ge 1}\cup\{\infty\}$,
$Q_\ell=[q_{\ell-1},q_\ell]\subset \Q_{\ge 1}\cup\{\infty\}$ is either a point or a closed interval,
$\mu_\ell:Q_\ell\to\Q\cup\{\infty\}$ is an affine function, non-constant when $Q_\ell$ is not a point, such that $\mu_\ell(q)\le q$ for all $q\in Q_\ell$ and $\min_{q\in Q_\ell}\mu_\ell(q)=\beta_\ell$.}
\end{definition}
\begin{definition}\label{equivalent_pizza}
\emph{Two pizzas are \emph{combinatorially equivalent} if the corresponding abstract pizzas are the same.}
\end{definition}

\begin{theorem}\label{contact_equiv} \emph{(See \cite[Theorem 4.9]{birbrair2014lipschitz}.)}
Two non-negative semialgebraic Lipschitz functions $f$ and $g$ defined on a H\"older triangle $T$ are contact Lipschitz equivalent
if, and only if, minimal pizzas on $T$ associated with $f$ and $g$ are combinatorially equivalent.
\end{theorem}

Let $T=T(\gamma_1,\gamma_2)$ and $T'=T(\gamma'_1,\gamma'_2)$ be Normally Embedded $\beta$-H\"older triangles satisfying the condition :
\begin{equation}\label{tord-tord}
tord(\gamma_1,T')=tord(\gamma_1,\gamma'_1)=tord(\gamma'_1,T),\quad tord(\gamma_2,T')=tord(\gamma_2,\gamma'_2)=tord(\gamma'_2,T).
\end{equation}

For example, the triangles $(T, Graph(f))$ considered in this section satisfy this condition. The following question is natural. Let $T=T(\gamma_1,\gamma_2)$ and $T'=T(\gamma'_1,\gamma'_2)$ be Normally Embedded semialgebraic $\beta$-H\"older triangles satisfying (\ref{tord-tord}).
Is it true that the union $T\cup T'$ is Lipschitz outer equivalent to the union $T\cup Graph(f)$, where $f$ is the distance function: $f(x)=dist(x,T')$? In the paper \cite{BG} the authors show that it is not true.

\section{Outer Lipschitz Geometry, Snakes}\label{section:snakes}

To formulate the results of this section, we need several definitions.

 \begin{definition}\label{DEF: Lipschitz non-singular arc}
	\normalfont Let $(X,0)$ be a surface germ. An arc $\gamma$ of $X$ is \textit{Lipschitz non-singular} if there
exists a Normally Embedded H\"older triangle $T \subset X$ such that $\gamma$ is an interior arc of  $T$ and $\gamma \not\subset
\overline{X\setminus T}$. Otherwise, $\gamma$ is \textit{Lipschitz singular}. It follows from the pancake decomposition that a
surface germ $X$ contains finitely many Lipschitz singular arcs.
The union of all Lipschitz singular arcs in $X$ is denoted by $Lsing(X)$.
A H\"older triangle $T\subset X$ is \textit{non-singular} if all interior arcs of $T$ are Lipschitz non-singular. 
\end{definition}

\begin{definition}\label{DEF: generic arcs}
\normalfont If $T=T(\gamma_1,\gamma_2)$ is a non-singular $\beta$-H\"older triangle, an arc $\gamma$ of $T$ is \textit{generic} if $itord(\gamma_1,\gamma) = itord(\gamma,\gamma_2) = \beta$. The set of generic arcs of $T$ is denoted $G(T)$. 
\end{definition}

\begin{definition}\label{DEF: abnormal arcs}
	\normalfont An arc $\gamma$ of a Lipschitz non-singular $\beta$-H\"older triangle $T$ is \textit{abnormal} if there are two normally embedded H\"older triangles $T'\subset T$ and $T''\subset T$ such that $T'\cap T'' = \gamma$ and $T\cup T'$ is not normally embedded.
Otherwise $\gamma$ is \textit{normal}. The set $Abn(T)$ of abnormal arcs of $T$ is outer Lipschitz invariant. 
\end{definition}

\begin{definition}\label{Def:snake}
	\normalfont A non-singular $\beta$-H\"older triangle $T$ is called a $\beta$\textit{-snake} if $Abn(T)=G(T)$. 
\end{definition}

The following important property of snakes can be interpreted as ``separation of scales'' in outer Lipschitz Geometry.

\begin{lemma}\label{Lem:minimal pancake decomp of a snake}
	Let $T$ be a $\beta$-snake, and let $\{T_k\}_{k=1}^p$ be a minimal pancake decomposition of $T$.
Then each $T_k$ is a $\beta$-H\"older triangle. 
\end{lemma}

\begin{figure}
	\centering
	\includegraphics[width=4.5in]{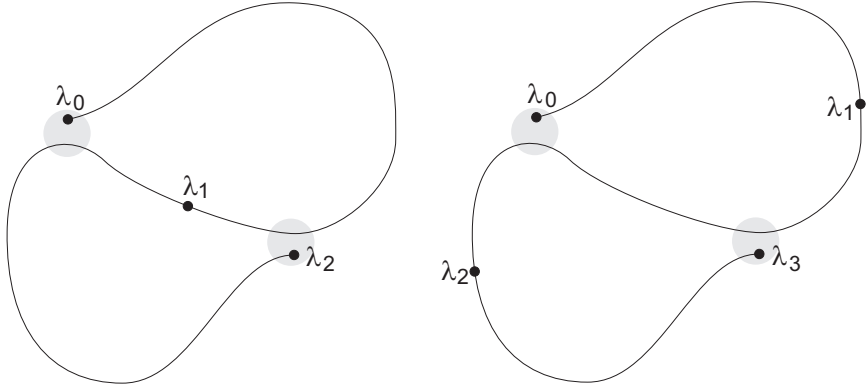}
	\caption{Two combinatorially non-equivalent minimal pancake decompositions of a snake. Black dots indicate the boundary arcs of pancakes.}\label{fig:Two pancake decompositions of a snake} 
\end{figure}

\begin{remark}\label{REM: non-equivalent minimal pancake decompositions}
\normalfont Minimal pancake decompositions of a snake may be combinatorially non-equivalent, as shown in Figure
\ref{fig:Two pancake decompositions of a snake}.
We use a planar plot to represent the link of a snake.
The points in Figure \ref{fig:Two pancake decompositions of a snake} correspond to arcs of the snake.
 The points with smaller Euclidean distance inside the shaded disks correspond to arcs with the tangency order higher than their inner tangency order $\beta$. Black dots indicate the boundary arcs of pancakes. 
\end{remark}

\begin{definition}\label{Def: weak LNE}
\normalfont A $\beta$-H\"older triangle $T$ is \textit{weakly normally embedded} if, for any two arcs $\gamma$ and $\gamma'$ of $T$ such that $tord(\gamma, \gamma') > itord(\gamma, \gamma')$, we have $itord(\gamma, \gamma')=\beta$. 
\end{definition}

\begin{proposition}\label{Prop:weak LNE}
	Let $T$ be a $\beta$-snake. Then $T$ is weakly normally embedded. 
\end{proposition}

Weak Lipschitz equivalence of snakes is a combinatorial invariant \cite[Subsection 6.3]{GS}.

\section{Tangent Cones}

\begin{definition}\label{tangent-cone}\normalfont
The {\it tangent cone} $C_0 X$ of a semialgebraic set $X$ at $0$ is defined as follows:
$$C_0 X=Cone\left(\,\lim_{\epsilon\to 0}\,\frac1\epsilon\,\big(X \cap\{|x|=\epsilon\}\big)\right),$$ where the limit here means the Hausdorff limit.

\end{definition}

\begin{remark}\label{hausdorff}\normalfont
There are several equivalent definitions of the tangent cone of a semialgebraic set. In particular, the tangent cone $C_0 X$ can be
defined as the set of tangent vectors at the origin to all the arcs in $X$. The tangent cone of a semialgebraic set is semialgebraic.
\end{remark}

The tangent cone is Lipschitz invariant:

\begin{theorem}\label{Sampaio}\emph{(See \cite{S}.)} If two germs $(X,0)$ and $(Y,0)$ are outer (resp. ambient)
Lipschitz equivalent, then the corresponding tangent cones $C_0 X$ and $C_0 Y$ are  outer (resp. ambient)
Lipschitz equivalent.
\end{theorem}

The result is used Theory of Metric Knots (see \cite{met-knots1},\cite{met-knots2}) to prove Universality Theorem below (see also
\cite{met-knots3},\cite{met-knots2}). This result was also used to prove that a complex analytic set, which is a Lipschitz submanifold of $\C^n$,y is a smooth submanifold (See\cite{regularity}). Moreover, the result was used in the recent study of the Zariski Multiplicity Conjecture (see the paper of Fernandes and Sampaio at the present volume).

\section{Ambient Equivalence. Metric Knots.}

\begin{definition}\label{equivalences}
\emph{Two germs of semialgebraic sets $(X,0)$ and $(Y,0)$ are \emph{outer Lipschitz equivalent}
if there exists a homeomorphism $H: (X,0)\rightarrow (Y,0)$ bi-Lipschitz with respect to the outer metric. The germs are \emph{
semialgebraic outer Lipschitz equivalent} if the map $H$ can be chosen to be semialgebraic.
The germs are \emph{ambient Lipschitz equivalent} if there exists an orientation preserving bi-Lipschitz homeomorphism
$\widetilde{H}: (\R^4,0) \rightarrow (\R^4,0)$, such that $\widetilde{H}(X)=Y$. The germs are \emph{ semialgebraic ambient
Lipschitz equivalent} if the map $\widetilde{H}$ can be chosen to be semialgebraic.}
\end{definition}
\begin{definition}\label{link}
\emph{The \emph{link at the origin} $L_X$ of a germ $X$ is the equivalence class of the sets $X\cap S^3_{0,\varepsilon}$ for small
positive $\varepsilon$ with respect to the ambient Lipschitz equivalence. The \emph{tangent link} of $X$ is the link at the origin
of the tangent cone of $X$.}
\end{definition}
\begin{remark}\label{well-defined}
\emph{By the finiteness theorems of Mostowski, Parusinski and Valette (see \cite{Mostowski}, \cite{Parusinski} and
\cite{valette2005Lip}) the link at the origin is well defined.
We write ``the link at the origin'' speaking of this notion of the link from Singularity Theory, reserving
the word ``link'' for the notion of the link in Knot Theory.
If $X$ has an isolated singularity at the origin then each connected component of $L_X$ is a knot in $S^3$.}
\end{remark}

The following result (so called {\bf Universality Theorem}) shows the difference between outer and ambient Lipschitz Geometry of
germs of real surfaces:

\begin{theorem}\emph{({\bf Universality Theorem.})}\label{universality}
Let $K \subset {S}^3$ be a knot. Then one can associate to $K$ a semialgebraic surface germ $(X_K,0)$ in $\R^4$ so that the
following holds:\newline
$1) \ $ The link at the origin of each germ $X_K$ is a trivial knot;\newline
$2) \ $ All germs $X_K$ are outer Lipschitz equivalent;\newline
$3) \ $ Two germs $X_{K_1}$ and $X_{K_2}$ are ambient semialgebraic Lipschitz equivalent only if the knots $K_1$ and $K_2$ are
isotopic.
\end{theorem}

To show how to proof works, we include the figure \ref{fig:xk}:

\begin{figure}
\centering
\includegraphics[width=5in]{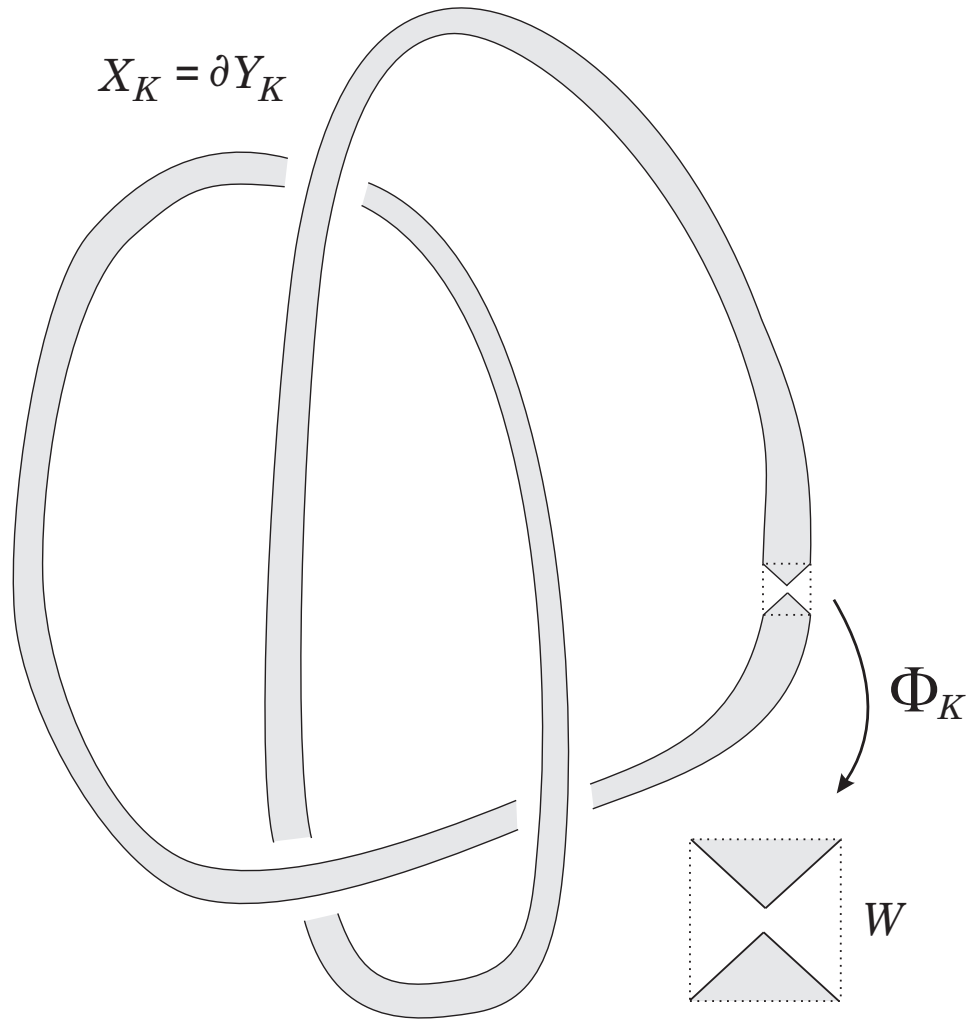}
\caption{The proof of Theorem \ref{universality}.}\label{fig:xk}
\end{figure}

The figure $4$, representing the link of the surface $X_K$. A reader can find a detailed explanation in \cite{met-knots3}.

 The following result shows that for given tangent cone one can find infinitely many Lipschitz outer equivalent, but not Lipschitz ambient equivalent, surface germs.

\begin{theorem}\label{twoknots-wedge}
For any two knots $K$ and $L$
there exists a semialgebraic surface germ $\tilde X_{KL}$ such that:

1. For any knots $K$ and $L$, the link at the origin of $\tilde X_{KL}$ is isotopic to $L$.

2. For any knots $K$ and $L$, the tangent link of $\tilde X_{KL}$ is isotopic to $K$.

3. For fixed $\alpha$ and $\beta$, all surface germs $\tilde X_{KL}$ are outer bi-Lipschitz equivalent.
\end{theorem}

\end{document}